\documentclass[twocolumn]{IEEEtran}

\IEEEoverridecommandlockouts    
\usepackage{psfrag}
\usepackage{amsmath,amssymb,amsfonts,bbm}
\usepackage{latexsym}
\usepackage{graphicx}

\usepackage{colortbl}
\usepackage{fancyhdr}
\usepackage{epstopdf}
\usepackage{float}
\usepackage{hyperref}
\usepackage{booktabs}
\usepackage{enumerate}
\usepackage[ruled]{algorithm2e}
\usepackage{balance}
\usepackage{mathrsfs}

\usepackage[usenames,dvipsnames,svgnames]{xcolor}

\newtheorem{theorem}{Theorem}
\newtheorem{definition}{Definition}
\newtheorem{proposition}{Proposition}

\newtheorem{corollary}{Corollary}

\newtheorem{remark}{Remark}

\newtheorem{standing}{Standing Assumption}

\IEEEoverridecommandlockouts

\newcommand{\R}{\mathbb{R}}

\newcommand{\mc}{\mathcal}

\newcommand{\col}{\mathrm{col}}

\newcommand{\dom}{\mathrm{dom}}
\newcommand{\proj}{\mathrm{proj}}

\newcommand{\bs}{\boldsymbol}

\graphicspath{{Figures/} }

\newcommand{\Rmnum}[1]{\expandafter\@slowromancap\romannumeral #1@}

\def\be{\begin{equation}}
\def\ee{\end{equation}}
\def\ba{\begin{array}}
\def\ea{\end{array}}
\def\eqa{\begin{eqnarray}}
\def\eqe{\end{eqnarray}}

\newtheorem{lem}{Lemma}

\usepackage[prependcaption,colorinlistoftodos]{todonotes}

\begin{document}

\title{Continuous-time integral dynamics for monotone aggregative games with coupling constraints}

\author{Claudio De Persis and Sergio Grammatico 
\thanks{C. De Persis is with the Faculty of Science and Engineering, University of Groningen, The Netherlands. S. Grammatico is with the Delft Center for Systems and Control (DCSC), TU Delft, The Netherlands. E-mail addresses: \texttt{c.de.persis@rug.nl}, \texttt{s.grammatico@tudelft.nl}. This work was partially supported by the Netherlands Organisation for Scientific Research (NWO) under research projects OMEGA (grant n. 613.001.702) and P2P-TALES (grant n. 647.003.003). A preliminary version of part of the results in this paper is in \cite{depersis:grammatico:18ecc}.
}
}

\thispagestyle{empty}
\pagestyle{empty}

\maketitle

\begin{abstract}
We consider continuous-time equilibrium seeking in monotone aggregative games with coupling constraints. We propose semi-decentralized integral dynamics and prove their global convergence to a variational generalized aggregative or Nash equilibrium. The proof is based on Lyapunov arguments and invariance techniques for differential inclusions.
\end{abstract}

{\keywords
Aggregative game theory, Multi-agent systems, Decentralized control, Projected dynamical systems.
}

\section{Introduction}

Aggregative game theory \cite{kukushkin:04} is a mathematical framework to model inter-dependent optimal decision making problems for a set of noncooperative agents, where the decision of each agent is affected by some aggregate effect of all the agents. Motivated by application domains where this aggregative feature arises, e.g.\ demand side management \cite{Saad2012} and network congestion control, equilibrium seeking in aggregative games is currently an active research area. 

Existence and uniqueness of (Nash) equilibria in (aggregative) games has been comprehensively studied, especially in close connection with variational inequalities \cite{facchinei:pang}, \cite[\S12]{palomar:eldar}. Distributed and semi-decentralized algorithms \cite{koshal:nedic:shanbhag:16,yi:pavel:17cdc}, \cite{grammatico:parise:colombino:lygeros:16, grammatico:17, belgioioso:grammatico:17cdc} have been proposed as \textit{discrete-time} dynamics that converge to an equilibrium of the game, e.g.\ Nash or aggregative equilibrium, under appropriate technical assumptions and sufficient conditions on the problem data. 
Specifically, one can characterize the desired equilibria as the zeros of a monotone operator, e.g.\ via the concatenation of interdependent Karush--Kuhn--Tucker operators, and formulate an equivalent fixed-point problem, to be solved via discrete-time dynamics with guaranteed global asymptotic convergence \cite{grammatico:17,belgioioso:grammatico:17cdc}.

Within the literature on equilibrium seeking for aggregative games with coupling constraints, the available solution methods are algorithms in \textit{discrete time}, where tuning the step size is typically a hard task. Therefore, in this paper, we address the aggregative equilibrium computation problem via \textit{continuous-time} dynamics. 
Inspired by passivity arguments \cite{pavel}, our original contribution is to provide simple, primal-dual integral dynamics for the computation of generalized aggregative and Nash equilibria via semi-decentralized dynamics.

To handle both local and global constraints, we propose equilibrium seeking dynamics that are characterized as the dynamics of a projected dynamical system \cite{nagurney:zhang}, whose solutions are intended as locally absolutely continuous functions.
Thus, we exploit invariance arguments for differential inclusions with maximally monotone set-valued right-hand side, and apply it to our primal-dual projected dynamics \cite{goebel.scl17, brogliato.scl}. From the technical perspective, our main contribution is to prove global asymptotic convergence of the proposed dynamics to a generalized (primal-dual) equilibrium of the aggregative game, under mild assumptions on the problem data, namely, local convexity of cost functions, convexity of constraints, and strict monotonicity of the pseudo-gradient mapping. Compared to our preliminary contribution \cite{depersis:grammatico:18ecc}, in this paper, we consider aggregative games with \textit{coupling constraints}, propose \textit{primal-dual} dynamics, and discuss convergence to both generalized aggregative equilibria and generalized Nash equilibria.

The paper is organized as follows. We introduce and mathematically characterize the problem setup in Section \ref{sec:aggregative.games}. We propose the equilibrium seeking dynamics and present the main result in Section \ref{sec:integral.dynamics}. Technical discussions and corollaries are in Section \ref{sec:discussion}. The proofs are given in the Appendix.

\subsubsection*{Notation and definitions}
$\bs{0}$ denotes a matrix/vector with all elements equal to $0$. $\otimes$ denotes the Kronecker product.
Given $N$ vectors $x_1, \ldots, x_N \in \R^n$, we define $\boldsymbol{x} := \col\left(x_1,\ldots,x_N\right) = \left[ x_1^\top, \ldots , x_N^\top \right]^\top$, $\boldsymbol{x}_{-i} := \col\left(x_1,\ldots , x_{i-1}, x_{i+1},\ldots,x_N\right)$, and $\textrm{avg}( \bs{x} ) := \frac{1}{N}\sum_{i=1}^{N} x_i$. 
Let the set $\mathcal{S} \subseteq \R^n$ be non-empty. 
The symbol ${\rm bdry}(\mathcal{S})$ denotes the boundary of $\mathcal{S}$, and 
the mapping $\iota_{\mc{S}}:\R^n \rightarrow \{ 0, \, \infty \}$ denotes the indicator function, i.e., $\iota_{\mc{S}}(x) = 0$ if $x \in \mathcal{S}$, $\infty$ otherwise. 
The set-valued mapping $\mathcal{N}_{\mc{S}}: \R^n \rightrightarrows \R^n$ denotes the normal cone operator, i.e., 
$\mathcal{N}_{\mc{S}}(x) = \varnothing$ if $x \notin \mc{S}$, $\left\{ v \in \R^n \mid \sup_{z \in \mc{S}} \, v^\top (z-x) \leq 0  \right\}$ otherwise. The set-valued mapping $\mathcal{T}_{\mc{S}}: \R^n \rightrightarrows \R^n$ denotes the tangent cone operator.
The mapping $\proj_{\mathcal{S}}(\cdot) := \textrm{argmin}_{ y \in \mathcal{S}} \left\| y - \cdot\right\| : \R^n \rightarrow \mathcal{S}$ denotes the projection operator; 
$\Pi_{ \mathcal{S} }(x,v) := \lim_{ h \rightarrow 0^{+} } \tfrac{1}{h}\left( \proj_{ \mathcal{S}}( x + h v ) - x  \right) $ denotes the projection of the vector $v \in \R^n$ onto the tangent cone of $\mathcal{S}$ at $x \in \mathcal{S}$, i.e., $\Pi_{\mathcal{S}}( x, \cdot ) = \proj_{ \mathcal{T}_{\mathcal{S}}(x) }(\cdot)$.
For a function $f: \R^n \rightarrow \overline{\R}$, $\dom(f) := \{x \in \R^n \mid f(x) < \infty\}$; $\partial f: \dom(f) \rightrightarrows {\R}^n$ denotes its subdifferential set-valued mapping, defined as $\partial f(x) := \{ v \in \R^n \mid f(z) \geq f(x) + v^\top (z-x)  \textup{ for all } z \in \textup{dom}(f) \}$; if $f$ is differentiable at $x$, then $\partial f(x) = \left\{ \nabla f(x)\right\}$. 
Given a closed convex set $C \subseteq \R^n$ and a single-valued mapping $F: C \rightarrow \R^n$, the variational inequality problem, denoted by VI$(C,F)$, is the problem to find $x^* \in C$ such that $ \inf_{ y \in C } \, (y-x^*)^\top  \, F(x^*) \geq 0.$

\section{Mathematical background: Aggregative games and variational equilibria}
\label{sec:aggregative.games}

\subsection{Jointly-convex aggregative games with coupling constraints}

A jointly-convex aggregative game with coupling constraints is denoted by a triplet $\mathcal{G}_{\textup{agg}} = \left( \mathcal{I}, ( J_i )_{i \in \mathcal{I}}, (\mathcal{X}_i)_{i \in \mathcal{I}} \right)$, where $\mathcal{I}:=\{1,\ldots, N\}$ is the index set of $N$ decision makers, or agents, 
$\left( J_i: \R^n \times \R^n \rightarrow \overline{\R} \right)_{i \in \mathcal{I}}$ is an ordered set of cost functions and $ \left( \mathcal{X}_i: \R^{ n(N-1) } \rightrightarrows \R^n \right)_{i \in \mathcal{I}}$ is an ordered set of set-valued mappings that represent coupled constraint sets. 
For each $i \in \mathcal{I}$, we assume an affine structure for the set $\mathcal{X}_i$:
$$\mathcal{X}_i( \bs{x}_{-i} ) := \{ y \in \Omega_i \mid A_i y + \textstyle \sum_{j \in \mathcal{I} \setminus \{i\}} A_j x_j \leq b \}\,,$$
for some set $\Omega_i \subseteq \R^n$ and matrices $A_1, \ldots, A_N \in \R^{ m \times n}$.

In aggregative games, the aim of each agent $i \in \mathcal{I}$ is to minimize its objective function $\bs{x} \mapsto J^i (x^i, {\rm avg}(\bs{x}))$ that depends on the local decision variable and on the average among the decision variables of all agents, i.e., ${\rm avg}( \bs{x} ) := \frac{1}{N}\sum_{i=1}^{N} x_i$. Formally, a jointly-convex aggregative game with coupling constraints represents the following collection of inter-dependent optimization problems:
\begin{equation}
\label{eq:Gagg}
\forall i \in \mathcal{I}: \, 
\left\{
\begin{array}{cl}
\displaystyle \min_{ x_i \in \Omega_i } & J_i\left( x_i, \textrm{avg}(\bs{x}) \right) \\
\textrm{ s.t.} & A \bs{x} - b \leq 0 \,,
\end{array}
\right.
\end{equation}
where $A \bs{x} := \left[ A_1, \ldots, A_N \right] \bs{x} = A_i x_i + \sum_{j \neq i} A_j x_j$.

Throughout the paper, we have the following assumption.

\smallskip
\begin{standing} \label{ass:regularity}
\textit{Continuity, compactness, convexity}.
The objective functions $\{J_i\}_{i \in \mathcal{I}}$ are continuous. The sets $\{\Omega_i\}_{i \in \mathcal{I}}$ are non-empty, compact and convex. The set $\bs{X} := \bs{\Omega} \cap \{ \bs{x} \in \R^{n N} \mid A \bs{x} \leq b \}$, where $\bs{\Omega} := \Omega_1 \times \ldots \times \Omega_N$, is non-empty and satisfies Slater's constraint qualification.
For all $i \in \mathcal{I}$, and $z \in \R^n$, the function $J_i( \cdot, z )$ is continuously differentiable and convex.
{\hfill $\square$}
\end{standing}

\subsection{Generalized aggregative equilibrium and strictly-monotone pseudo-gradient mapping}

Our aim is to design continuous-time, semi-decentralized, dynamics that asymptotically converge to a generalized aggregative equilibrium, which is a set of decision variables such that each is optimal given the average among all the decision variables and the coupling constraints.

\smallskip
\begin{definition} \label{def:GAE}
\textit{Generalized aggregative equilibrium}.
A set of decision variables $\bs{x}^* = {\rm col}\left({x}_1^{*} , \ldots , {x}_{N}^* \right) \in \bs{X}$ is a generalized aggregative equilibrium of the game in \eqref{eq:Gagg} if, for all $i \in \mathcal{I}$,  
$$ 
\begin{array}{rcl}
x_i^* \in & \displaystyle \underset{y \in \Omega_i}{ \textrm{argmin} } &  
J_i\left( y , \textrm{avg}(\bs{x}^*) \right) \\
 & \textup{ s.t. } & A_i y + \sum_{j \neq i} A_j x_j^* \leq b.
\end{array}
\vspace{-0.5cm}
$$
{\hfill $\square$}
\end{definition}
\smallskip

A fundamental mapping in game theory is the so-called pseudo-subdifferential mapping, which in our setup with continuously differentiable functions, hence single-valued subdifferentials, is a pseudo-gradient mapping. Since we are interested in generalized aggregative equilibria, rather than generalized Nash equilibria, together with semi-decentralized equilibrium seeking dynamics, let us adopt the following definition of (semi-extended) pseudo-gradient mapping: 
\begin{equation}\label{eq:F}
F( \bs{x}, \sigma ) := \left[
\begin{matrix}
\mathrm{col}\left(\left( \nabla_{x_i} J_i( x_i, \sigma ) \right)_{i \in \mathcal{I}}\right) \\
k \left( \sigma - \textrm{avg}( \bs{x} ) \right)
\end{matrix}
\right]\,,
\end{equation}
where $k > 0$ is a free design parameter, and $\sigma$ is a control variable. Throughout the paper, we assume that the pseudo-gradient mapping in \eqref{eq:F} is strictly monotone, see the discussion in Section \ref{sec:discussion-smon} for sufficient conditions on the problem data.

\smallskip
\begin{standing} \label{ass:strictly-monotone-F}
\textit{Strictly-monotone pseudo-gradient mapping}.
The pseudo-gradient mapping $F$ in \eqref{eq:F} is strictly monotone on $\bs{X} \times \R^n$, i.e., for all $\xi,\zeta \in \bs{X} \times \R^n$ such that $\xi \neq \zeta$, $\langle F(\xi) - F(\zeta) , \xi-\zeta \rangle > 0.$
{\hfill $\square$}
\end{standing}
\smallskip

It follows by the regularity of the problem data (Standing Assumption \ref{ass:regularity}) and by the strict monotonicity of the pseudo-gradient mapping (Standing Assumption \ref{ass:strictly-monotone-F}) that a generalized aggregative equilibrium exists. The proof is analogous to that of \cite[Th. 2.3.3 (a)]{facchinei:pang}.

\subsection{Operator-theoretic characterization}

With the aim to decouple the coupling constraints of the game, $A \bs{x} \leq b$ in \eqref{eq:Gagg}, we adopt duality theory for equilibrium problems. We start from the definition of the Lagrangian functions, $\{ L_i \}_{i \in \mathcal{I}}$, one for each agent $i \in \mathcal{I}$:
\begin{equation}\label{eq:Lagrangian}
L_i\left( x_i, \bs{x}_{-i}, \lambda_i \right) := J_i( x_i, \textrm{avg}( \bs{x} ) ) + {\lambda_i}^\top ( A \bs{x} - b ) \,,
\end{equation}
where $\lambda_i$ is a dual variable. Then, for each $i \in \mathcal{I}$, we introduce the Karush--Kuhn--Tucker (KKT) system:

\begin{equation}\label{eq:KKT}
\forall i \in \mathcal{I}: \ \left\{
\begin{array}{l}
0 = \nabla_{x_i} J_i( x_i, \sigma ) + \iota_{\Omega_i}(x_i) + A_i^\top \lambda_i \\
0 = \sigma - \textrm{avg}( \bs{x} ) \\
0 \leq \lambda_i \perp - (A \bs{x}-b) \geq 0 \,,
\end{array}
\right.
\end{equation}
where $\{\lambda_i\}_{i \in \mathcal{I}}$ are the dual variables, one vector for each agent $i \in \mathcal{I}$, associated with the coupling constraint, and $0 \leq \lambda_i \perp - (A \bs{x}-b) \geq 0$ represents the complementarity condition. Note that in \eqref{eq:KKT}, the first two equations are equivalent to $0 = \nabla_{x_i} J_i( x_i, \textrm{avg}( \bs{x} ) ) + \iota_{\Omega_i}(x_i) + A_i^\top \lambda_i$. We use the former formulation to recover a semi-decentralized solution algorithm later on.
Next, we follow the steps in \cite{belgioioso:grammatico:17cdc} and focus on {the class of \textit{variational}} generalized aggregative equilibria, i.e., generalized aggregative equilibria that satisfy the KKT system in \eqref{eq:KKT} with equal dual variables, $\lambda_i = \lambda$ for all $i \in \mathcal{I}$. 
{In our definition of variational generalized aggregative equilibrium, the connection with the solutions to the KKT system is inspired by \cite[Th. 9, Def. 3]{facchinei:kanzow:07}.} We refer to \cite[\S 5]{facchinei:kanzow:07} for the relevant properties of variational (generalized Nash) equilibria.

\smallskip
\begin{definition}\label{def:v-GAE}
\textit{Variational generalized aggregative equilibrium}.
A set of decision variables $\bs{x}^* = {\rm col}\left({x}_1^{*} , \ldots , {x}_{N}^* \right) \in \bs{X}$ is a variational generalized aggregative equilibrium (v-GAE) of the game in \eqref{eq:Gagg} if it is a GAE of the game in \eqref{eq:Gagg} and there exists $\lambda^* \in \R_{\geq 0}^m$ such that the triplet $\left( \bs{x}^*, \textrm{avg}(\bs{x}^*), \bs{1}_N \otimes \lambda^* \right)$ solves the KKT system in \eqref{eq:KKT}.
{\hfill $\square$}
\end{definition}
\smallskip

Existence and uniqueness of the v-GAE follows by the standing assumptions, due to the connection with variational inequalities -- the proof is analogous to \cite[Prop. 12.11]{palomar:eldar}. By introducing the dual variable $\lambda$, we have extended the space of the decision variables of the aggregative game. It then follows that the extended version of the pseudo-gradient mapping, 
\begin{equation}\label{eq:extended-pseudo-gradient}
F_{\textup{ext}}( \bs{x}, \sigma, \lambda ) := 
\left[
\begin{matrix}
F(\bs{x},\sigma) \smallskip  \\
b
\end{matrix}
\right] + 
\left[
\begin{matrix}
\phantom{-}\bs{0} & \bs{0} & A^\top \\
\phantom{-}\bs{0} & \bs{0} & \bs{0} \\
-A     & \bs{0} & \bs{0}
\end{matrix}
\right] 
\left[ 
\begin{matrix}
\bs{x} \\
\sigma \\
\lambda
\end{matrix}
\right] \,,
\end{equation}
has a fundamental role in the operator-theoretic characterization of the equilibria. Specifically, we show in the following that the solution of the KKT system is a zero of a (maximally) monotone operator that contains the extended pseudo-gradient mapping in \eqref{eq:extended-pseudo-gradient} and that it generates a v-GAE.

\smallskip
\begin{lem} \label{lem:agg-eq=zero} 
\textit{Operator-theoretic characterization}. The following statements are equivalent:
\begin{enumerate}[(i)]
\item $\bs{x}^*$ is a v-GAE of the game in \eqref{eq:Gagg}; 
\item $\left( \bs{x}^*, \textrm{avg}( \bs{x}^* ), \lambda^* \right) \in \textrm{zer}\left( \mathcal{N}_{ \bs{\Omega} \times \R^n \times \R_{\geq 0}^m } + F_{\textup{ext}}\right)$, 
for some $\lambda^* \in \R_{\geq 0}^{m}$.
\hfill $\square$
\end{enumerate}
\end{lem}
\smallskip

\begin{remark}\label{rem:max.monotonicity} \textit{Maximal monotonicity}.
The normal cone $\mathcal{N}_{ \bs{\Omega} \times \R^n \times \R_{\geq 0}^m }$ is maximally monotone \cite[Ex. 20.26]{bauschke:combettes}, $F_{\textup{ext}}$ is monotone and continuous, hence maximally monotone \cite[Cor. 20.28]{bauschke:combettes}, and $\textrm{dom}\left( F_{\textup{ext}} \right) = \R^{nN} \times \R^n \times \R^m$. Thus, the sum operator $\mathcal{N}_{ \bs{\Omega} \times \R^n \times \R_{\geq 0}^m } + F_{\textup{ext}}$ in Lemma \ref{lem:agg-eq=zero} is maximally monotone as well \cite[Cor. 25.5 (i)]{bauschke:combettes}.
\hfill $\square$
\end{remark}

\section{Continuous-time integral dynamics for generalized aggregative equilibrium seeking}
\label{sec:integral.dynamics}

For asymptotically reaching the v-GNE, we consider the following continuous-time integral dynamics:
\begin{equation}
\label{aggregative.dynamics}
\ba{rccl}
\forall i \in \mathcal{I}: & \dot x_i &=& \Pi_{\Omega_i} \left( x_i \, , \, -\nabla_{x_i} J_i( x_i, \sigma ) - A_i^\top \lambda \right) \\
 & \dot \sigma  &=& k \left( {\rm avg}( \bs{x} ) -\sigma  \right) \\
 & \dot \lambda &=& \Pi_{\R_{\geq 0}^m}\left( \lambda ,  A \bs{x} - b \right).
\ea
\end{equation}
where $k > 0$ is a free parameter gain.

Equivalently, in collective projected-vector form, the dynamics in \eqref{aggregative.dynamics} read as
\begin{equation}
\label{aggregative.dynamics.vector}
\begin{bmatrix} 
\dot{\bs{x}} \\ 
\dot \sigma \\ 
\dot \lambda \end{bmatrix} = 
\Pi_{ \bs{\Omega} \times \R^n \times \R_{\geq 0}^m } 
\left( 
\left[\begin{matrix}
\bs{x} \\
\sigma \\ 
\lambda 
\end{matrix}\right]
 \, , \, \begin{bmatrix}
 - F( \bs{x}, \sigma ) +
 \left[
\begin{matrix}
 -A^\top \lambda  \\
 \bs{0}
\end{matrix} 
 \right]
 \\
 A \bs{x} - b
\end{bmatrix}  \right).
\end{equation}

\smallskip
\begin{remark}
\textit{Semi-decentralized structure}. 
The computation and information exchange in \eqref{aggregative.dynamics} are semi-decentralized: each agent performs decentralized computations, namely, projected-pseudo-gradient steps, and does not exchange information with other agents. A central control unit, which does not participate in the game, collects aggregative information, ${\rm avg}(\bs{x}(t))$ and $A \bs{x}(t)-b$, and broadcasts two signals, $\sigma(t)$ and $\lambda(t)$, to the agents playing the aggregative game. In turn, the dynamics of the broadcast signal $\sigma(t)$ are driven by the average among all the decision variables, ${\rm avg}(\bs{x}(t))$, while the dynamics of the signal $\lambda(t)$ are driven by the coupling-constraint violation, $A \bs{x}(t) - b$. The semi-decentralized structure prevents that the agents are imposed to exchange truthful information.
\hfill $\square$
\end{remark}
\smallskip

First, we show that the $\bs{x}-$part of an equilibrium for the dynamics in \eqref{aggregative.dynamics.vector} is a v-GAE, in view of Lemma \ref{lem:agg-eq=zero}.

\smallskip
\begin{lem} \label{lem:agg.eq.and.dyn.eq}
The following statements are equivalent:
\begin{enumerate}[(i)]
\item $\left(\bs{\bar x}, \bar \sigma, \bar{\lambda} \right)$ is an equilibrium for the dynamics in \eqref{aggregative.dynamics.vector};
\item $\left(\bs{\bar x}, \bar \sigma, \bar{\lambda} \right) \in \textrm{zer}\left( \mathcal{N}_{ \bs{\Omega} \times \R^n \times \R_{\geq 0}^m } + F_{\textup{ext}} \right)$.
\hfill $\square$
\end{enumerate}
\end{lem}
\smallskip

In view of Lemma \ref{lem:agg.eq.and.dyn.eq}, we can directly analyze the convergence of the projected dynamics in \eqref{aggregative.dynamics.vector} to an equilibrium. Let us introduce a quadratic  function, $V$, which is used later on to obtain a Lyapunov function.

\smallskip
\begin{lem} \label{prop:new-3}
Consider the function
\begin{multline} \label{eq:V}
V(\bs{x},\sigma,\lambda) := 
\tfrac{1}{2} \left\| \bs{x} \!-\! \bs{x}' \right\|^2 + \tfrac{1}{2} \left\| \sigma \!-\! \sigma'\right\|^2 + \tfrac{1}{2} \left\| \lambda \!-\! \lambda' \right\|^2 \,,
\end{multline}
where $(\bs{x},\sigma,\lambda), (\bs{x}',\sigma',\lambda')$ are arbitrary vectors in ${\bs{\Omega}} \times \mathbb{R}^n \times \R_{\geq 0}^m$.  It holds that
\begin{equation} \label{strict.lyap}
\ba{rl}
\dot V(\bs{x},\sigma,\lambda) := &\nabla V(\bs{x},\sigma,\lambda)^\top 
\left[ 
\begin{smallmatrix}
\dot{\bs{x}} \\
\dot{\sigma} \\
\dot{\lambda}
\end{smallmatrix}
\right] \\[3mm]
\le& 
\nabla V(\bs{x},\sigma,\lambda)^\top 
\begin{bmatrix}
 - F( \bs{x}, \sigma ) +
 \left[
\begin{matrix}
 -A^\top \lambda  \\
 \bs{0}
\end{matrix} 
 \right]
 \\
 A \bs{x} - b
\end{bmatrix},
\ea
\end{equation}
where $\left[ 
\begin{smallmatrix}
\dot{\bs{x}} \\
\dot{\sigma} \\
\dot{\lambda}
\end{smallmatrix}
\right]$ stands for the right-hand side in \eqref{aggregative.dynamics.vector}.
\hfill $\square$
\end{lem}
\smallskip

We are now ready to establish our main global asymptotic convergence result. The proof, given in Appendix A, is based on invariance arguments for differential inclusions with maximal monotone set-valued right-hand side.

\smallskip

\begin{theorem} \label{th:convergence-GAE}
\textit{Global asymptotic convergence to variational generalized aggregative equilibrium}. 
Let $\bs{x}^*$ be the v-GAE of the game in \eqref{eq:Gagg}. 
For any initial condition $\left( \bs{x}_0, \sigma_0, \lambda_0 \right) \in \bs{\Omega} \times \mathbb{R}^n \times \R_{\geq 0}^m$, there exists a unique solution to  \eqref{aggregative.dynamics.vector} starting from $\left( \bs{x}_0, \sigma_0, \lambda_0 \right)$, which is a locally absolutely continuous function satisfying \eqref{aggregative.dynamics.vector} almost everywhere,   remains in $ \bs{\Omega} \times \mathbb{R}^n \times \R_{\geq 0}^m$, is bounded for all time, and converges to
$\{ \bs{x}^* \} \times \{ {\rm avg}\left( \bs{x}^* \right) \} \times \{  \overline{\lambda}  \}$,
a Lyapunov stable equilibrium of \eqref{aggregative.dynamics.vector}.
\hfill $\square$
\end{theorem}

\section{Technical discussions} \label{sec:discussion}

\subsection{On the strict monotonicity of the pseudo-gradient} \label{sec:discussion-smon}

In this subsection, let us consider a separable structure for the cost functions, i.e.,
\begin{equation} \label{eq:Ji-separable}
\forall i \in \mathcal{I}: \ J_i(x_i, \sigma ) = f_i(x_i) + \left( C_i \sigma \right)^\top x_i \,,
\end{equation}
for some convex functions $\{f_i\}_{i\in \mathcal{I}}$ and $n \times n$ matrices $\{C_i\}_{i\in \mathcal{I}}$. 
The next result provides a sufficient condition on the problem data such that the pseudo-gradient mapping $F$ in \eqref{eq:F} is strictly monotone (Standing Assumption \ref{ass:strictly-monotone-F}).

\smallskip
\begin{proposition}\label{prop:suff.cond.strict.monoton}
Consider the aggregative game in \eqref{eq:Gagg} with cost functions $\{J_i\}_{i \in \mathcal{I}}$ as in \eqref{eq:Ji-separable}. Assume that, for each $i \in \mathcal{I}$, the function $f_i$ is twice continuously differentiable and $\mu-$strongly convex in $\Omega_i$. The pseudo-gradient mapping in \eqref{eq:F} is strictly monotone if 
\begin{equation} \label{eq:suff-cond}
\min\{\mu,k\} > \tfrac{1}{2}\left( \max_{i \in \mathcal{I}} \left\| C_i\right\|_{\infty} + \tfrac{k}{N}\right).
\vspace{-0.25cm}
\end{equation}
\hfill $\square$
\end{proposition}
\smallskip

\begin{proof}
By \eqref{eq:Ji-separable}, the pseudo-gradient mapping in \eqref{eq:F} is 
$$
F( \bs{x}, \sigma) = 
\mathrm{col}\left( \left( \partial f_i(x_i) + C_i \sigma \right)_{i=1}^{N} , k ( \sigma - \mathrm{avg}(\bs{x}) )  \right)\,,
$$
hence its subdifferential reads as
$$
\partial F( \bs{x},\sigma ) = 
\left[
\begin{array}{ccc|c}
\partial^2 f_1(x_1) &  &  & C_1 \\
 & \ddots & & \vdots \\
 & & \partial^2 f_N(x_N) & C_N \\ 
 \midrule
 -\tfrac{k}{N} I & \cdots &  -\tfrac{k}{N} I & k I
\end{array}
\right].
$$
It follows from \cite[Prop. 2.3.2]{facchinei:pang} that the pseudo-gradient $F$ is strictly monotone if and only if its subdifferential is positive semi-definite, i.e., since $f_i$ is $\mu-$strongly convex,
$$
\left[
\begin{array}{ccc|c}
\mu I &  &  & \tfrac{1}{2}\left( C_1 -\tfrac{k}{N} I\right) \\
 & \ddots & & \vdots \\
 & & \mu I & \tfrac{1}{2}\left( C_N -\tfrac{k}{N} I\right) \\
 \midrule
* & \cdots & * & k I
\end{array}
\right] \succ 0.
$$

By the Gershgorin circle theorem, the latter is true if $\min\{ \mu , k \} > \tfrac{1}{2} |c_{i,i}-\tfrac{k}{N}| + \tfrac{1}{2} \sum_{j \neq i} | c_{i,j} | $, which is implied by \eqref{eq:suff-cond}.
\end{proof}
\smallskip

The sufficient condition in \eqref{eq:suff-cond} extends that in \cite[Prop. 1]{depersis:grammatico:18ecc} to the case of heterogeneous matrices $\{C_i\}_{i\in \mathcal{I}}$. In turn, it is less restrictive than the sufficient condition in \cite[Th. 2]{grammatico:16cdc-convex}. The inequality condition in \eqref{strict.lyap} becomes less restrictive as $N$ grows, which is desirable for large number of agents \cite[\S IV]{depersis:grammatico:18ecc}.

\vspace{-0.5cm}

\subsection{On separable convex coupling constraints}

Let us discuss the setup with separable, nonlinear yet convex, coupling constraints, i.e., of the form
\begin{equation} \label{eq:separable-convex-coupling-constraints}
\textstyle g(\bs{x}) := \sum_{i=1}^{N} g_i( x_i ) \leq 0 \,,
\end{equation}
where the functions $\{ g_i \}_{i\in \mathcal{I}}$ are convex and continuously differentiable, and the set $\{ \bs{x} \mid g(\bs{x}) \leq 0 \}$ is non-empty and satisfies Slater's constraint qualification. To recover affine coupling constraints, the optimization problems of the agents can be rewritten with auxiliary decision variables as
\begin{equation}
\label{eq:Gagg-2}
\forall i \in \mathcal{I}: \, 
\left\{
\begin{array}{cl}
\displaystyle \min_{ x_i , y_i } & J_i\left( x_i, \textrm{avg}(\bs{x}) \right) \\
\textrm{ s.t.} & (x_i,y_i) \in \Omega_i \times \mathcal{Y}_i \\
               & g_i(x_i) \leq y_i \\
               & \sum_{j=1}^N y_j \leq 0 \,,
\end{array}
\right.
\end{equation}
where the set $\mathcal{Y}_i \supset \{ g_i(\xi) \in \R^m \mid \xi \in \Omega_i \} $ is compact and convex. Now, if $\bs{x}^*$ is the GAE of the original game with coupling constraints as in \eqref{eq:separable-convex-coupling-constraints}, then the pair $\left( \bs{x}^*, \bs{y}^* \right)$, with $y_i^{*} = g_i(x_i^*)$ for all $i$, is a GAE of the game in \eqref{eq:Gagg-2}. 
Conversely, let $\left( \bs{x}^*, \bs{y}^* \right)$ be a GAE of the game in \eqref{eq:Gagg-2}. If the coupling constraint is inactive at the equilibrium, $\sum_{j=1}^N y_j^* < 0$, then it is unnecessary and $\bs{x}^*$ is a GAE of the original game; if the coupling constraint is active, then
\begin{equation*}
\begin{array}{r}
\forall i \in \mathcal{I}: \   (x_i^*, y_i^*) \ \in   \\
  \\  
  \\
\end{array}
\left\{
\begin{array}{cl}
\displaystyle \underset{x_i,y_i}{\mathrm{argmin}} & J_i\left( x_i, \textrm{avg}(\bs{x}) \right) \\
\textrm{ s.t.} & (x_i,y_i) \in \Omega_i \times \mathcal{Y}_i \\
              & g_i(x_i) + \sum_{j \neq i} y_j^* \leq 0. \\ 
\end{array}
\right.
\end{equation*}
Therefore, the pair $\left( \bs{x}^*, ( g_i( x_i^* )  )_{i \in \mathcal{I}} \right)$ is a GAE of the game in \eqref{eq:Gagg-2}, and in turn $\bs{x}^*$ is a GAE of the original game.

\subsection{On generalized Nash equilibria}

We recall that a Nash equilibrium is a set of strategies where each is optimal given the other strategies, as formalized next.

\smallskip
\begin{definition} \label{def:GNE}
\textit{Generalized Nash equilibrium}.
A set of decision variables $\bs{x}^* = {\rm col}\left({x}_1^{*} , \ldots , {x}_{N}^* \right) \in \bs{X}$ is a generalized Nash equilibrium (GNE) of the game in \eqref{eq:Gagg} if, for all $i \in \mathcal{I}$,  
$$ 
\begin{array}{rcl}
x_i^* \in & \displaystyle \underset{y \in \Omega_i}{ \textrm{argmin} } &  
J_i\left( y , \tfrac{1}{N} y + \tfrac{1}{N} \sum_{j \neq i} x_j^*  \right) \\
 & \textup{ s.t. } & A_i y + \sum_{j \neq i} A_j x_j^* \leq b.
\end{array}
\vspace{-0.5cm}
$$
{\hfill $\square$}
\end{definition}
\smallskip

\begin{remark}
A GNE in Definition \ref{def:GNE} differs from a GAE in Definition \ref{def:GAE}, since in the latter, each decision variable is optimal given the average among the decision variables of all agents that enters as second argument of the cost functions. We refer to \cite{li:zhang:zhao:lian:kalsi:16cdc,deori:margellos:prandini:17} for a comparison between aggregative/mean-field equilibria and Nash equilibria.
{\hfill $\square$}
\end{remark}
\smallskip

If we aim at computing a GNE, rather than a GAE, then the definition of pseudo-gradient mapping shall be changed to
\begin{equation}\label{eq:F-GNE}
F_{\textrm{N}}( \bs{x}, \sigma ) := 
\left[
\begin{matrix}
\mathrm{col}\left( \left( \nabla_{x_i} J_1( x_i, \sigma ) + \tfrac{1}{N} \nabla_{\sigma} J_i( x_i, \sigma ) \right)_{i=1}^{N} \right) \\
k \left( \sigma - \textrm{avg}( \bs{x} )\right)
\end{matrix}
\right] \,,
\end{equation}
since, for each agent $i$, the variable $x_i$ enters as local decision variable in both the first and the second argument of the cost function $J_i$. Analogously to \eqref{aggregative.dynamics.vector}, possible continuous-time generalized Nash equilibrium seeking dynamics are

\begin{equation}
\label{aggregative.dynamics.Nash}
\begin{bmatrix} 
\dot{\bs{x}} \\ 
\dot \sigma \\ 
\dot \lambda \end{bmatrix} = 
\Pi_{ \bs{\Omega} \times \R^n \times \R_{\geq 0}^m } 
\left( 
\left[\begin{matrix}
\bs{x} \\
\sigma \\ 
\lambda 
\end{matrix}\right]
 \, , \, \begin{bmatrix}
 - F_{\textrm{N}}( \bs{x}, \sigma ) +
 \left[
\begin{matrix}
 -A^\top \lambda  \\
 \bs{0}
\end{matrix} 
 \right]
 \\
 A \bs{x} - b
\end{bmatrix}  \right).
\end{equation}

Convergence to a variational GNE (v-GNE) of the above dynamics then follows if the pseudo-gradient mapping is strictly monotone.

\smallskip
\begin{corollary} \label{cor:convergence-GNE}
\textit{Global asymptotic convergence to generalized Nash equilibrium}.
Let $\bs{x}_{\textrm{N}}^*$ be the v-GNE of the game in \eqref{eq:Gagg}.
Assume that the mapping $F_{\textrm{N}}$ in \eqref{aggregative.dynamics.Nash} is strictly monotone on $\bs{\Omega} \times \mathbb{R}^n \times \R_{\geq 0}^m$. For any initial condition $\left( \bs{x}_0, \sigma_0, \lambda_0 \right) \in \bs{\Omega} \times \mathbb{R}^n \times \R_{\geq 0}^m$, there exists a unique solution to  \eqref{aggregative.dynamics.Nash} starting from $\left( \bs{x}_0, \sigma_0, \lambda_0 \right)$, which is a locally absolutely continuous function satisfying \eqref{aggregative.dynamics.vector} almost everywhere,  remains in $ \bs{\Omega} \times \mathbb{R}^n \times \R_{\geq 0}^m$, is bounded for all time, 
and converges to $\{ \bs{x}_{\textrm{N}}^* \} \times \{{\rm avg}\left( \bs{x}_{\textrm{N}}^* \right)\} \times \{ \overline{\lambda} \}$, a Lyapunov stable equilibrium of \eqref{aggregative.dynamics.Nash}. 
\hfill $\square$
\end{corollary}
\smallskip

Analogously to Proposition \ref{prop:suff.cond.strict.monoton}, in the case of separable cost functions as in \eqref{eq:Ji-separable}, we provide sufficient conditions on the problem data such that the pseudo-gradient mapping $F_{\textrm{N}}$ in \eqref{aggregative.dynamics.Nash} is strictly monotone.

\smallskip
\begin{proposition}\label{prop:suff.cond.strictly.monotone.Nash}
Consider the aggregative game in \eqref{eq:Gagg} with cost functions $\{J_i\}_{i \in \mathcal{I}}$ as in \eqref{eq:Ji-separable}. Assume that, for each $i \in \mathcal{I}$, the function $f_i$ is twice continuously differentiable and $\mu-$strongly convex in $\Omega_i$. The mapping $F_{\textrm{N}}$ in \eqref{aggregative.dynamics.Nash} is strictly monotone if 
\begin{equation} \label{eq:suff-cond-Nash}
\min\{\mu,k\} > \max_{i \in \mathcal{I}} \, \tfrac{1}{2} \left\| C_i\right\|_{\infty} + \tfrac{1}{N}\left\| C_i \right\|_1 + \tfrac{1}{2} \tfrac{k}{N}.
\vspace{-0.25cm}
\end{equation}
\hfill $\square$
\end{proposition}
\smallskip

\begin{proof}
Since $\nabla_{x_i} J_i(x_i, \sigma) + \tfrac{1}{N} \nabla_{\sigma} J_i(x_i, \sigma) = \partial f_i(x_i) + C_i \sigma + \tfrac{1}{N} C_i^\top x_i $, by the proof of Proposition \ref{prop:suff.cond.strict.monoton}, we shall have 
$$
\left[
\begin{array}{ccc|c}
\mu I + \tfrac{1}{N} C_1^{\top} &  &  & \tfrac{1}{2}\left( C_1 -\tfrac{k}{N} I\right) \\
 & \ddots & & \vdots \\
 & & \mu I + \tfrac{1}{N} C_N^{\top} & \tfrac{1}{2}\left( C_N -\tfrac{k}{N} I\right) \\
 \midrule
* & \cdots & * & k I
\end{array}
\right] \succ 0.
$$

Thus, by the Gershgorin circle theorem, the latter is true if $\min\{ \mu + \tfrac{1}{N} c_{i,i} , k \} > \tfrac{1}{N} \sum_{j \neq i} | c_{j,i} | + \tfrac{1}{2} |c_{i,i}-\tfrac{k}{N}| + \tfrac{1}{2} \sum_{j \neq i} | c_{i,j} | $, which is implied by \eqref{eq:suff-cond-Nash}.
\end{proof}

\section{Conclusion and outlook} \label{sec:conclusion}

In aggregative games with affine coupling constraints, continuous-time integral dynamics with semi-decentralized computation and information exchange can ensure global asymptotic convergence to generalized aggregative or Nash equilibria, under mild regularity and strict monotonicity assumptions. Future research will focus on continuous-time distributed-averaged-integral dynamics in multi-agent network games with coupling constraints.

\section*{Appendix A: Proofs}

For ease of notation, next, we use $\xi := \textrm{col}(\bs{x},\sigma)$, $\xi^* := \textrm{col}(\bs{x}^*,\sigma^*)$, $\bar \xi := \textrm{col}(\bar {\bs{x}},\bar \sigma)$  and $\Xi := \bs{\Omega} \times \R^n$.

\smallskip
\subsubsection*{Proof of Lemma \ref{lem:agg-eq=zero}} 
Analogous to \cite[Prop. 12.4]{palomar:eldar}, namely, with $F$ in \eqref{eq:F}, rather than $F_{\textrm{N}}$ in \eqref{eq:F-GNE}, see \cite[Th. 1]{belgioioso:grammatico:17cdc}.
\hfill $\blacksquare$

\smallskip
\subsubsection*{Proof of Lemma \ref{lem:agg.eq.and.dyn.eq}}
By Moreau's decomposition theorem,
\begin{multline*}
\bs{0} = \Pi_{ \Xi }\left(  
\bar{\xi}
 , \, -F( \bar{\xi}  )  +  
\left[ 
\begin{smallmatrix}
- A^\top \bar{\lambda} \\
\bs{0}
\end{smallmatrix} 
\right] 
  \right) \\
  = -F(\bar{\xi} , \bar{\sigma} ) +  
\left[ 
\begin{smallmatrix}
- A^\top \bar{\lambda} \\
\bs{0}
\end{smallmatrix} 
\right] \, -\proj_{ \mathcal{N}_{ \bs{\Omega} \times \R^n }( \bar{\xi}  )}\left( -F( \bar{\xi}  ) +  
\left[ 
\begin{smallmatrix}
- A^\top \bar{\lambda} \\
\bs{0}
\end{smallmatrix} 
\right] \right)
\end{multline*}
and
$\bs{0} = \Pi_{ \R_{\geq 0}^m }\left( \bar{\lambda}, A \bar{\bs{x}} - b \right) = A \bs{x} - b - \proj_{ \mathcal{N}_{ \R_{\geq 0}^n }( \bar\lambda ) }\left( A \bar{\bs{x}} - b\right)$.
The proof then follows immediately.
\hfill $\blacksquare$

\subsection*{Proof of Lemma \ref{prop:new-3}}
The proof follows the steps of \cite[Proof of Lemma 6]{pavel}.
Since ${\nabla V}(\xi,\lambda)^\top = \mathrm{col}\left( \xi - \xi', \lambda - \lambda' \right)^\top$, 
for all vectors $\bs{u}$, by Moreau's decomposition theorem, we have that
\begin{multline*}
\left( \xi - \xi'\right)^\top \Pi_{ \Xi }\left( \xi, -F(\xi) + \bs{u} \right) = \\
\left( \xi - \xi'\right)^\top \left[ -F(\xi) + \bs{u} \,  - 
\textrm{proj}_{ \mathcal{N}_{\Xi }(\xi)}( -F(\xi) + \bs{u} ) \right].
\end{multline*} 
By definition of the normal cone $\mathcal{N}_{\Xi}(\xi)$, we have that
\begin{equation*}
-\left( \xi - \xi' \right)^\top \textrm{proj}_{ \mathcal{N}_{\Xi }(\xi)}( -F(\xi) + \bs{u} ) \leq 0,
\end{equation*}
and in turn
\begin{equation}\label{eq:diss-in1}
\left( \xi - \xi'\right)^\top \Pi_{ \Xi }\left( \xi, -F(\xi) + \bs{u} \right) \leq \left( \xi - \xi'\right)^\top \left( -F(\xi) + \bs{u} \right).
\end{equation}

With similar arguments, we can show that 
\begin{equation}\label{eq:diss-in2}
\left( \lambda - \lambda'\right)^\top \Pi_{ \mathbb{R}^m_{\ge 0} }\left( \lambda, A\bs{x}-b\right) \leq \left( \lambda - \lambda'\right)^\top  \left( A\bs{x}-b \right).
\end{equation}

The proof follows by summing up the inequalities in \eqref{eq:diss-in1} with $\bs{u} = \left[ 
\begin{smallmatrix}
-A^\top \lambda \\
\bs{0}
\end{smallmatrix}
\right]$ and \eqref{eq:diss-in2}.
\hfill $\blacksquare$

\subsection*{Proof of Theorem \ref{th:convergence-GAE}}
The dynamics in \eqref{aggregative.dynamics.vector} represent a projected dynamical system with discontinuous right-hand side \cite{nagurney}. The proof uses  invariance arguments for differential inclusions with maximally monotone right-hand side \cite{goebel.scl17}. First, we note that $F_{\textup{ext}}$ in \eqref{eq:extended-pseudo-gradient} is continuous and monotone. Then, we consider a zero of $\mathcal{N}_{ 
\Xi \times \R_{\geq 0}^m}+F_{{\rm ext}}$ (Lemma \ref{lem:agg-eq=zero}), $\left( \xi^*, \lambda^* \right)$, and, bearing in mind Lemma \ref{prop:new-3}, define the Lyapunov function $W( \xi, \lambda ) := \tfrac{1}{2}\left\| \xi - \xi^*\right\|^2 + \tfrac{1}{2}\left\| \lambda - \lambda^*\right\|^2$. We show next that 
\begin{equation}
\label{aux-vi}
\nabla W(z)^\top F_{\textup{ext}}(z^*)=
\begin{bmatrix}
\xi-\xi^*\\
\lambda-\lambda^*
\end{bmatrix}^\top
F_{\textup{ext}}( \xi^*, \lambda^* ) 
\ge 0 
\end{equation}
for all $z = (\xi, \lambda) \in \Xi \times \R_{\geq 0}^m$. By Lemma \ref{lem:agg-eq=zero} and \ref{lem:agg.eq.and.dyn.eq}, 
\begin{align}
\label{eq:equilib.identity} 
\begin{split}
\bs{0} &= \Pi_{ \Xi }\left( \left[ 
\begin{matrix}
{\bs{x}}^* \\
 \sigma^*
\end{matrix}
 \right] , \, -F( {\bs{x}}^*, {\sigma}^* ) +  
\left[ 
\begin{matrix}
- A^\top {\lambda}^* \\
\bs{0}
\end{matrix} 
\right] 
  \right) \\
  \bs{0} &= \Pi_{ \R_{\geq 0}^m }\left( {\lambda}^*, A {\bs{x}}^* - b \right)\,,
  \end{split}
\end{align}
therefore, we have 
$0 = -\nabla W(\xi, \lambda)^\top \left[ 
\begin{smallmatrix}
\dot \xi \\
\dot{\lambda}
\end{smallmatrix}
\right]^*$, 
where $\left[ 
\begin{smallmatrix}
\dot \xi\\
\dot{\lambda}
\end{smallmatrix}
\right]^* $ 
stands for the right-hand side of \eqref{eq:equilib.identity}. By Lemma \ref{prop:new-3}, we immediately obtain \eqref{aux-vi}:
$$
0 = -\nabla W(\xi, \lambda)^\top \left[ 
\begin{smallmatrix}
\dot \xi\\
\dot{\lambda}
\end{smallmatrix}
\right]^*
\le 
\nabla W(\xi,\lambda)^\top 
F_{\textup{ext}}(\xi^*, \lambda^* ).
$$

Consequently, we have that
\[\ba{rl}
\nabla W(z)^\top \dot z = & -\nabla W(z)^\top F_{\textup{ext}}(z)\\
\le & 
-\nabla W(z)^\top [F_{\textup{ext}}(z)-F_{\textup{ext}}(z^* )]\le 0,
\ea\]
by the monotonicity of $F_{\textup{ext}}$. We conclude that $W$ is not increasing along the trajectories of  \eqref{aggregative.dynamics.vector}. By radial unboundedness of  $W$, for any initial condition $z_0$, the corresponding solution is bounded and therefore the associated $\omega$-limit set $\Lambda(z_0)$ is non-empty, compact, invariant and attractive. Moreover, by definition of the $\omega$-limit set, $W$ is constant on $\Lambda(z_0)$. Thus, any solution $\zeta(\cdot)$ with initial condition in $\Lambda(z_0)$ must satisfy $\dot W( \zeta(t))=0$, that is $\Lambda(z_0)$ is contained in the set of points satisfying $\nabla W(z)^\top F_{\textup{ext}}(z) = 0$. 
We then study the set $\mathcal{O} = \left\{ z \in  \Xi \times \R_{\geq 0}^m  \mid  
\nabla W(z)^\top F_{\textup{ext}}(z) = 0 \right\}$.
For all $\col(\xi,\lambda) \in \mathcal{O}$, it holds:
\begin{multline}
\nabla W(\xi, \lambda)^\top F_{\textup{ext}}( \xi,\lambda ) = 
\left[
\begin{matrix}
\xi - \xi^* \\
\lambda - \lambda^*
\end{matrix}
\right]^\top F_{\textup{ext}}( \xi,\lambda ) \\
= \left( \xi - \xi^* \right)^\top\left( F(\xi) + 
\left[
\begin{matrix}
A^\top \lambda \\
\bs{0}
\end{matrix}
\right]
 \right) - 
 \left( \lambda - \lambda^* \right)^\top \left( A \bs{x} - b\right).
\end{multline}

By Lemma \ref{lem:agg-eq=zero}, we have that 
$F_{\textup{ext}}( \xi^*,\lambda^* ) + 
\left[ 
\begin{smallmatrix}
v^* \\
\bs{0} \\
\end{smallmatrix}
\right] = \bs{0}
$
for some $v^* \in \mathcal{N}_{\bs{\Omega}}(\bs{x}^*)$, hence 
$F(\xi^*) + 
\left[
\begin{smallmatrix}
A^\top \lambda^* \\
\bs{0}
\end{smallmatrix}
\right] + 
\left[
\begin{smallmatrix}
v^* \\
\bs{0}
\end{smallmatrix}
\right] = \bs{0}$ 
and ${\lambda^*}^\top ( A \bs{x}^* - b ) = 0$. Therefore, for all $\col(\xi,\lambda) \in \mathcal{O}$,
\begin{multline}\label{eq:equality.-1}
0 =  \left( \xi - \xi^* \right)^\top \left( F(\bs{x}, \sigma) -F(\bs{x}^*, \sigma^*) 
- 
\left[
\begin{matrix}
v^* \\
\bs{0}
\end{matrix}
\right]
+ \right. \\
\left.
\left[ 
\begin{matrix}
A^\top ( \lambda - \lambda^* )  \\
\bs{0}
\end{matrix}
\right]
\right) - \left( \lambda - \lambda^* \right)^\top \left( A \bs{x} - b \right).
\end{multline}

Now, we observe that 
$\left( \lambda \!-\! \lambda^* \right)^\top \left( A \bs{x}^* \!-\! b \right) = 
\underbrace{\lambda^\top}_{ \geq 0 } \underbrace{\left( A \bs{x}^* \!-\! b \right)}_{ \leq 0 } - \underbrace{{\lambda^*}^\top \left( A \bs{x}^* \!-\! b \right)}_{ = 0 } \leq 0 \,,$
and in turn
\begin{align}
\label{eq:equality.0}
\begin{split}
0 \geq & \left[ 
\begin{matrix}
\bs{x} - \bs{x}^* \\
\sigma - \sigma^*
\end{matrix}
\right]^\top \left( F(\bs{x}, \sigma) -F(\bs{x}^*, \sigma^*) - 
\left[
\begin{matrix}
v^* \\
\bs{0}
\end{matrix}
\right] + \right. \\
 & \left.
\left[
\begin{matrix}
A^\top ( \lambda - \lambda^* )  \\
\bs{0}
\end{matrix} 
\right] \right) - \left( \lambda - \lambda^* \right)^\top A \left( \bs{x}- \bs{x}^*  \right) \\
 = &  \left( \xi - \xi^* \right)^\top \left( F(\bs{x}, \sigma) -F(\bs{x}^*, \sigma^*) - 
 \left[
\begin{matrix}
v^* \\
\bs{0}
\end{matrix}
\right]
 \right) \, \geq 0.
 \end{split}
\end{align}
The last inequality holds because, by Standing Assumption \ref{ass:strictly-monotone-F}, $\left(  F(\bs{x}, \sigma) -F(\bs{x}^*, \sigma^*) \right)^\top \left( \xi - \xi^* \right) \geq 0$ and, by the definition of normal cone, ${v^*}^\top (\bs{x} - \bs{x}^*) \leq 0$. Thus, we obtain 
\begin{equation}
\label{eq:equality.1}
\left(  F(\bs{x}, \sigma) -F(\bs{x}^*, \sigma^*) \right)^\top \left( \xi - \xi^* \right) =  {v^*}^\top (\bs{x} - \bs{x}^*) = 0.
\end{equation}

From \eqref{eq:equality.1}, due to Standing Assumption \ref{ass:strictly-monotone-F}, we conclude that $\bs{x} = \bs{x}^*$ and $\sigma=\sigma^*={\rm avg}(\bs{x}^*)$.
From \eqref{eq:equality.-1} and \eqref{eq:equality.1}, we obtain 
$0=\left( \bs{x}- \bs{x}^*  \right)^\top A^\top \left( \lambda - \lambda^* \right)- \left( \lambda - \lambda^* \right)^\top \left( A \bs{x} - b \right)$, 
hence $\lambda^\top \left( A \bs{x}^* - b \right) = 0$. 
The latter implies $(\lambda'-\lambda)^\top \left( A \bs{x}^* - b \right) \le 0$ for all $\lambda'\in \mathbb{R}^m_{\ge 0}$, i.e., $A \bs{x}^* - b\in \mathcal{N}_{\mathbb{R}^m_{\ge 0}}(\lambda)$, or, equivalently, $\bs{0}=\Pi_{\mathbb{R}^m_{\ge 0}}(\lambda, A \bs{x}^* - b)$.  The latter and the identity $\xi=\xi^*$ established before returns that $(\xi^*, \lambda)$ is a zero of $\mathcal{N}_{ \Xi \times \R_{\geq 0}^m} +F_{{\rm ext}}$, hence an equilibrium of \eqref{aggregative.dynamics.vector}, and this concludes the characterisation of $\mathcal{O}$.

We finally show that convergence is to an equilibrium point of \eqref{aggregative.dynamics.vector}. By Lemma \ref{lem:th1-brogliato} in Appendix B, the solution to  \eqref{aggregative.dynamics.vector} is the same as the solution to $-\dot z \in  F_{{\rm ext}}(z)+\mathcal{N}_{\Xi\times\mathbb{R}^m_{\ge 0}}(z)$, where the right-hand side of the differential inclusion is maximally monotone by Remark \ref{rem:max.monotonicity}. We can then apply \cite[Ch. 3, Sec. 2, Th. 1]{aubin-cellina}, \cite[Th. 2.2, (C1), (C3)]{goebel.scl17}, to conclude that every equilibrium of \eqref{aggregative.dynamics.vector} is Lyapunov stable and that, if the solution has an $\omega$-limit point at an equilibrium, then the solution converges to that equilibrium. Now, from the arguments in the first part of the proof, the non-empty and invariant $\omega$-limit set $\Lambda({\rm col}\left( \xi_0, \lambda_0 \right))$ is contained in $\mathcal{O}$. Since points of $\mathcal{O}$ are equilibria of \eqref{aggregative.dynamics.vector}, then the $\omega$-limit set $\Lambda({\rm col}\left( \xi_0, \lambda_0 \right))$ is a singleton with an  equilibrium $(\xi^*, \lambda)$ to which the solution converges.  This concludes the proof. 
\hfill $\blacksquare$
\smallskip

\smallskip
\textit{Proof of Corollary \ref{cor:convergence-GNE}}: 
Analogous to the proof of Theorem \ref{th:convergence-GAE}, namely, with $F_{\textrm{N}}$ in \eqref{eq:F-GNE}, in place of $F$ in \eqref{eq:F}.
\hfill $\blacksquare$

\section*{Appendix B: Projected dynamical systems}\label{app:B}
We consider a generic projected dynamical system 
\be\label{eq:pds}
\ba{rl}
\dot{z}= & \Pi_{K} \left(z , 
 - F( z ) \right)
\ea
\ee
where $K \subseteq \mathbb{R}^n$ is a non-empty, closed and convex set. 
The dynamic behavior of \eqref{eq:pds} is well-studied for continuous, hypomonotone mappings $F$.

\smallskip
\begin{definition}\label{asspt:hypomonotone}
\textit{Hypomonotonicity}.
A mapping $F:\mathbb{R}^n\to \mathbb{R}^n$ is hypomonotone if there exists $\beta \geq 0$ such that
\begin{equation*}
(z-z')^\top (F(z)-F(z'))\ge - \beta \|z-z'\|^2 
\end{equation*}
for all $z, z'\in \mathbb{R}^n$. 
\hfill $\square$
\end{definition}
\smallskip

In view of \cite{brogliato.scl}, \cite{nagurney}, we recall next some equivalent formulations of the projected dynamical system in \eqref{eq:pds}.

\smallskip
\begin{lem}[{from \cite[Th. 1]{brogliato.scl}}]
\label{lem:th1-brogliato}
Let $F$ in \eqref{eq:pds} be continuous and hypomonotone. 
For any initial condition $z_0\in K$, the differential inclusion
\begin{equation}
\label{udi.bis}
-\dot{z}(t) \stackrel{\textup{a.e.}}{\in} F(z(t)) + \mathcal{N}_K(z(t)). 
\end{equation}
has a unique solution $z(t)$ that belongs to $K$ for almost all $t\ge 0$.
Furthermore, the evolution variational inequality 
\begin{equation} \label{evi.bis}
z(t) \stackrel{\textup{a.e.}}{\in} K, \ t \ge 0 , \, 
\underset{ v \in K }{\inf} \langle \dot{z}(t)+ F(z(t)), v- z(t)\rangle \stackrel{\textup{a.e.}}{\ge} 0 \,,
\end{equation}
and the projected dynamical system 
\begin{equation*}
\dot {z}(t) \stackrel{\textup{a.e.}}{=} \mathrm{proj}_{\mathcal{T}_K(z(t))} (-F(z(t))) =  \Pi_K( z(t), -F(z(t)) )
\end{equation*}
have the same solution as to \eqref{udi.bis}. 
\hfill $\square$
\end{lem}

\balance

\bibliographystyle{IEEEtran}
\bibliography{library,int-nask-seeking}

\begin{thebibliography}{10}
\providecommand{\url}[1]{#1}
\csname url@samestyle\endcsname
\providecommand{\newblock}{\relax}
\providecommand{\bibinfo}[2]{#2}
\providecommand{\BIBentrySTDinterwordspacing}{\spaceskip=0pt\relax}
\providecommand{\BIBentryALTinterwordstretchfactor}{4}
\providecommand{\BIBentryALTinterwordspacing}{\spaceskip=\fontdimen2\font plus
\BIBentryALTinterwordstretchfactor\fontdimen3\font minus
  \fontdimen4\font\relax}
\providecommand{\BIBforeignlanguage}[2]{{%
\expandafter\ifx\csname l@#1\endcsname\relax
\typeout{** WARNING: IEEEtran.bst: No hyphenation pattern has been}%
\typeout{** loaded for the language `#1'. Using the pattern for}%
\typeout{** the default language instead.}%
\else
\language=\csname l@#1\endcsname
\fi
#2}}
\providecommand{\BIBdecl}{\relax}
\BIBdecl

\bibitem{depersis:grammatico:18ecc}
{C. De Persis} and S.~Grammatico, ``Continuous-time integral dynamics for
  aggregative game equilibrium seeking,'' in \emph{Proc. of the {IEEE} European
  Control Conference}, 2018.

\bibitem{kukushkin:04}
N.~S. Kukushkin, ``Best response dynamics in finite games with additive
  aggregation,'' \emph{Games and Economic Behavior}, vol.~48, no.~1, pp.
  94--10, 2004.

\bibitem{Saad2012}
W.~Saad, Z.~Han, H.~Poor, and T.~Ba\c{s}ar, ``Game theoretic methods for the
  smart grid,'' \emph{IEEE Signal Processing Magazine}, pp. 86--105, 2012.

\bibitem{facchinei:pang}
F.~Facchinei and J.~Pang, \emph{Finite-dimensional variational inequalities and
  complementarity problems}.\hskip 1em plus 0.5em minus 0.4em\relax Springer
  Verlag, 2003.

\bibitem{palomar:eldar}
D.~Palomar and Y.~Eldar, \emph{Convex optimization in signal processing and
  communication}.\hskip 1em plus 0.5em minus 0.4em\relax Cambridge University
  Press, 2010.

\bibitem{koshal:nedic:shanbhag:16}
J.~Koshal, A.~Nedi\'c, and U.~Shanbhag, ``Distributed algorithms for
  aggregative games on graphs,'' \emph{Operations Research}, vol.~64, no.~3,
  pp. 680--704, 2016.

\bibitem{yi:pavel:17cdc}
P.~Yi and L.~Pavel, ``A distributed primal-dual algorithm for computation of
  generalized {Nash} equilibria via operator splitting methods,'' \emph{Proc.
  of the IEEE Conf. on Decision and Control}, pp. 3841--3846, 2017.

\bibitem{grammatico:parise:colombino:lygeros:16}
S.~Grammatico, F.~Parise, M.~Colombino, and J.~Lygeros, ``Decentralized
  convergence to {Nash} equilibria in constrained deterministic mean field
  control,'' \emph{IEEE Trans. on Automatic Control}, vol.~61, no.~11, pp.
  3315--3329, 2016.

\bibitem{grammatico:17}
S.~Grammatico, ``Dynamic control of agents playing aggregative games with
  coupling constraints,'' \emph{IEEE Trans. on Automatic Control}, vol.~62,
  no.~9, pp. 4537 -- 4548, 2017.

\bibitem{belgioioso:grammatico:17cdc}
G.~Belgioioso and S.~Grammatico, ``Semi-decentralized {Nash} equilibrium
  seeking in aggregative games with coupling constraints and non-differentiable
  cost functions,'' \emph{IEEE Control Systems Letters}, vol.~1, no.~2, pp.
  400--405, 2017.

\bibitem{pavel}
D.~Gadjov and L.~Pavel, ``A passivity-based approach to {N}ash equilibrium
  seeking over networks,'' 2017, arXiv:1705.02424.

\bibitem{nagurney:zhang}
A.~Nagurney and D.~Zhang, \emph{Projected dynamical systems and variational
  inequalities with applications}.\hskip 1em plus 0.5em minus 0.4em\relax
  Springer, 1996.

\bibitem{goebel.scl17}
R.~Goebel, ``Stability and robustness for saddle-point dynamics through
  monotone mappings,'' \emph{Systems \& Control Letters}, vol. 108, pp. 16--22,
  2017.

\bibitem{brogliato.scl}
B.~Brogliato, A.~Daniilidis, C.~Lemarechal, and V.~Acary, ``On the equivalence
  between complementarity systems, projected systems and differential
  inclusions,'' \emph{Systems \& Control Letters}, vol.~55, no.~1, pp. 45--51,
  2006.

\bibitem{facchinei:kanzow:07}
F.~Facchinei and C.~Kanzow, ``Generalized {Nash} equilibrium problems,''
  \emph{A Quarterly Journal of Operations Research}, vol.~5, pp. 173--210,
  2007.

\bibitem{bauschke:combettes}
H.~H. Bauschke and P.~L. Combettes, \emph{Convex analysis and monotone operator
  theory in {Hilbert} spaces}.\hskip 1em plus 0.5em minus 0.4em\relax Springer,
  2010.

\bibitem{grammatico:16cdc-convex}
S.~Grammatico, ``Aggregative control of large populations of noncooperative
  agents,'' in \emph{Proc. of the IEEE Conf. on Decision and Control}, Las
  Vegas, USA, 2016.

\bibitem{li:zhang:zhao:lian:kalsi:16cdc}
S.~Li, W.~Zhang, L.~Zhao, J.~Lian, and K.~Kalsi, ``On social optima of
  non-cooperative mean field games,'' in \emph{Proc. of the IEEE Conf. on
  Decision and Control}, 2016, pp. 3584--3590.

\bibitem{deori:margellos:prandini:17}
L.~Deori, K.~Margellos, and M.~Prandini, ``On the connection between {Nash}
  equilibria and social optima in electric vehicle charging control games,'' in
  \emph{Proc. of the IFAC World Congress}, Toulouse, France, 2017.

\bibitem{nagurney}
A.~Nagurney and D.~Zhang, \emph{Projected dynamical systems and variational
  inequalities with applications}.\hskip 1em plus 0.5em minus 0.4em\relax
  Kluwer Academic Publishers, 1996.

\bibitem{aubin-cellina}
J.-P. Aubin and A.~Cellina, \emph{Differential inclusions}.\hskip 1em plus
  0.5em minus 0.4em\relax Springer-Verlag, 1984.

\end{thebibliography}

\end{document}